\documentclass[twoside,11pt]{article}
\usepackage{rotating}
\usepackage[all,cmtip]{xy}
\usepackage{amsfonts,amssymb,amsmath}
\usepackage[dvips]{epsfig}
\usepackage{booktabs}
\newcommand{\affiliation}[1]{\let\thefootnote\relax\footnote{\mbox{}\\ \noindent {#1}}}

\begin{document}
\title{\bf The Structure Constants of the Exceptional Lie Algebra ${\mathfrak g}_2$ in the Cartan-Weyl Basis}
\author{H. Fakhri\thanks{Email: hfakhri@tabrizu.ac.ir}\,, M.
Sayyah-Fard\thanks{Email: msayyahfard@tabrizu.ac.ir}\,, S.
Laheghi\thanks{Email: slaheghi@tabrizu.ac.ir}
\\
{\small {\em Department of Theoretical Physics and
Astrophysics, Faculty of Physics,}}\\
{\small {\em University of Tabriz, P. O. Box 51666-16471, Tabriz,
Iran}}} \maketitle
\begin{abstract}\noindent The purpose of this paper is to answer the
question whether it is possible to realize simultaneously the
relations $N_{\alpha,\beta}=-N_{-\alpha,-\beta}$,
$N_{\alpha,\beta}=N_{\beta,-\alpha-\beta}=N_{-\alpha-\beta,\alpha}$
and
$N_{\alpha,\beta}N_{-\alpha,-\beta}=-\frac{1}{2}q(p+1)\langle\alpha,H_{\alpha}\rangle$
by the structure constants of the Lie algebra ${\mathfrak g}_2$. We
show that if the structure constants obey the first relation, the
three last ones are violated, and vice versa. Contrary to the second
case, the first one uses the Cartan matrix elements to derive the
structure constants in the form of $\langle\beta,H_{\alpha}\rangle$.
The commutation relations corresponding to the first case are
exactly documented in the prior literature. However, as expected, a
Lie algebra isomorphism is established between the Cartan-Weyl bases
obtained in both approaches.
\\ \\{\it Mathematics Subject Classification
(2010): 17B22; 22E60}
\\ \\{\it Physics and Astronomy Classification Scheme (2008):
02.20.Sv}\\ \\
{\it Keywords:} Exceptional Lie algebra ${\mathfrak g}_2$;
Semisimple Lie algebras; Cartan-Weyl basis; Cartan-Chevalley basis;
Root system; Killing form; Commutation relations; Structure
constants
\end{abstract}

\section{Introduction}
The simple Lie algebras over the complex numbers have been
classified into four infinite classes, including the special linear,
odd orthogonal, symplectic and even orthogonal algebras, as well as
that into five finite exceptional Lie algebras ${\mathfrak g}_2$,
${\mathfrak f}_4$, ${\mathfrak e}_6$, ${\mathfrak e}_7$ and
${\mathfrak e}_8$ with the dimensions 14, 56, 78, 133 and 248,
respectively
\cite{Killing1,Killing2,Killing3,Killing4,Cartan1,Cartan2,Cartan3}.
Among the exceptional Lie algebras, the smallest and easiest one to
consider the root system, the Dynkin diagram, fundamental
representations, the structure constants, etc. is the algebra
${\mathfrak g}_2$. Its corresponding Lie group as a subgroup of the
Spinor group $Spin(7)$ is obtained by fixing a point in $S^7$
\cite{Adams}. The group $G_2$ has been attracted great attentions
in many areas such as gauge theories in particle physics
\cite{Loginov,Wipf1,Wipf2,Maas,Saearp,Shnir}, (super)gravity theory
\cite{House1,House2,Herfray} and etc. For example, in order to add a
new matter to the minimal $SU(3)$ model for electroweak unification,
it is developed by embedding in the exceptional Lie group $G_2$
\cite{Carone}. The reason to choose this is that any such group must
be at least rank 2 and contain the subgroup $SU(3)$. $G_2$ is not
only a rank 2 group that contains $SU(3)$ but also it has the
minimum number of generators with respect to any other group that
properly involves $SU(3)$.

Many studies have investigated different aspects of the Lie group
$G_2$ and the Lie algebra ${\mathfrak g}_2$. Some areas are listed
here. The complex form of $G_2$ is the isotropy group of a generic
3-form in 7-dimensional complex space \cite{Bincer1,Grigorian}.
Whilst, the compact and noncompact real forms of the group $G_2$
appear as the automorphism groups of the octonion algebra and the
split octonions, respectively \cite{Cartan2,Cartan3}. In Refs.
\cite{Jacobson1,Jacobson2}, it has been shown that the generators of
${\mathfrak g}_2$ can be appeared in the role of the derivations of
a certain non-associative Cayley algebra. The Casimir operators of
${\mathfrak g}_2$, with the generators written in terms of $A_2$ and
$B_3$ bases, have been considered in Refs. \cite{Bincer1,Bincer2}.
All finite-dimensional irreducible representations of ${\mathfrak
g}_2$ have been realized in spaces of complex-valued polynomials of
six variables \cite{Saenkarun1}.

Our goal here is to derive the commutation relations corresponding
to the exceptional Lie algebra ${\mathfrak g}_2$ in the Cartan-Weyl
basis from two different approaches, one is regularly used and the
other is unknown. We focus on the relations that should be satisfied
by the structure constants corresponding to the generators of the
type of raising and lowering operators. We denote that the structure
constants of the two bunches of commutation relations saturate the
known equations but different from each other. However, it will be
shown, by using an appropriate set of the linear relations between
the two bunches of the Cartan-Weyl bases, that the commutation
relations of the two different viewpoints are isomorphic to each
other. Let us follow the problem in details in the next section.

\section{Preliminaries}
Let ${\mathfrak g}$ be a semisimple, complex and finite-dimensional
Lie algebra. Suppose ${\mathfrak h}$ is a maximal abelian subalgebra
of ${\mathfrak g}$ and $\Delta$ is the system of all non-zero roots
of ${\mathfrak g}$ with respect to it. The non trivial subalgebra
${\mathfrak h}$ is called the Cartan subalgebra and consists of
linearly independent and simultaneously $\mbox{ad}$-diagonalizable
elements. This yields the triangular decomposition of ${\mathfrak
g}$ with respect to ${\mathfrak h}$ \cite{Onishchik,Kerf,Knapp},
\begin{eqnarray}
{\mathfrak g}={\mathfrak
h}\oplus\left(\bigoplus_{\alpha\in\Delta}{\mathfrak
g}_{\alpha}\right),
\end{eqnarray}
in which ${\mathfrak g}_{\alpha}$'s are eigenspaces of the operators
belonging to $\mbox{ad} {\mathfrak h}$,
\begin{eqnarray}\label{adhe}
\mbox{ad} H (E_{\alpha})\equiv[H,E_{\alpha}]=\langle\alpha,H\rangle
E_{\alpha}, \hspace{10mm}H\in{\mathfrak h},
\hspace{1mm}E_{\alpha}\in{\mathfrak g}_{\alpha}.
\end{eqnarray}
The Cartan subalgebra ${\mathfrak h}$ is an eigenspace with the
eigenvalue $0$ and so ${\mathfrak h}={\mathfrak g}_{0}$. The
eigenvalue equation (\ref{adhe}) offers a definition for the root
$\alpha$ as the linear form $\alpha(H)\equiv\langle\alpha,H\rangle$
on ${\mathfrak h}$. The main properties of roots and root subspaces
are listed as following:
\begin{itemize}
\item If $\alpha$  is a root belonging to $\Delta$, then so is
$-\alpha$. For $-1\neq a\in\mathbb{C}$, $a\alpha$ does not belong to
the root system $\Delta$. This allows us to partition root system
$\Delta$ into two subsets $\Delta_+$ and $\Delta_-$, containing the
positive and negative roots, respectively:
$\Delta=\Delta_+\bigcup\Delta_-$.
\item For $\alpha,\beta,\alpha+\beta\in\Delta$ we have
$[{\mathfrak g}_{\alpha},{\mathfrak g}_{\beta}]\subset{\mathfrak
g}_{\alpha+\beta}$ and $[{\mathfrak g}_{\alpha},{\mathfrak
g}_{-\alpha}]\subset{\mathfrak h}$. Also, $[{\mathfrak
g}_{\alpha},{\mathfrak g}_{\beta}]=0$ if $\alpha+\beta$ is not a
root of $\Delta$. \item All root subspaces ${\mathfrak g}_{\alpha}$
are one-dimensional. Therefore, we can assign a basis $E_{\alpha}$
to any root subspaces ${\mathfrak g}_{\alpha}$.
\item If $\alpha+\beta\neq0$, then $E_{\alpha}$ and $E_{\beta}$ are orthogonal with respect
to the Killing form
$K(E_{\alpha},E_{\beta})\equiv\mbox{Tr}(\mbox{ad}E_{\alpha}\circ\mbox{ad}E_{\beta})$.
\item The restriction of $K$ to ${\mathfrak h}\times{\mathfrak h}$ is
non-degenerate. For any $\alpha$, there is a unique element
$H_{\alpha}\in{\mathfrak h}$ such that
$K(H_{\alpha},H)=\langle\alpha,H\rangle$ for all $H\in{\mathfrak
h}$. \item The Killing form is nondegenerate on ${\mathfrak
g}\times{\mathfrak g}$. The invariance of the Killing form fixes the
normalization of the $E_{\pm\alpha}$ generators to one,
$K(E_{\alpha},E_{-\alpha})=1$. Then, we have
$[E_{\alpha},E_{-\alpha}]=H_{\alpha}$.
\end{itemize}
From the above considerations for ${\mathfrak g}$, the commutation
relations in the Cartan-Weyl basis are resulted as
\begin{eqnarray}\label{HHEH}
&&\hspace{-10mm}[H_{\alpha},H_{\beta}]=0, \hspace{10mm}
[H_{\alpha},E_{\beta}]=\langle\beta,H_{\alpha}\rangle
E_{\beta},\,\,\,\,\alpha,\beta\in\Delta,\\\label{EEH}
&&\hspace{-10mm}[E_{\alpha},E_{-\alpha}]=H_{\alpha}, \hspace{5mm}
[E_{\alpha},E_{\beta}]=N_{\alpha,\beta}E_{\alpha+\beta},\,\,\,\,\,\alpha+\beta\neq0,\,\,\,\,\,\,\alpha,\beta\in\Delta,
\end{eqnarray}
in which, $N_{\alpha,\beta}=0$ if $\alpha+\beta$ is not a root. The
root $\alpha\in\Delta$ is called simple if it can not be written as
a sum of two positive roots. Let us collect the properties of simple
roots \cite{Vilenkin}:
\begin{itemize}
\item Simple roots are linearly independent and their number
coincides with the rank of ${\mathfrak g}$, which is, the dimension
$l$ of the Cartan subalgebra ${\mathfrak h}$. We denote the set of
simple roots by $\pi=\{\alpha_1,\alpha_2, \cdots,
\alpha_l\}\subset\Delta_+$. \item The difference of two simple roots
is not a root. \item Every positive root is a linear combination of
the simple roots with non-negative integer coefficients.
\end{itemize}
For every simple root $\alpha_i$ one can consider the standard
triple $\{H_{\alpha_i}, E_{\alpha_i}, F_{\alpha_i}\}$ in ${\mathfrak
g}$ with $H_{\alpha_i}\in{\mathfrak h}$, $E_{\alpha_i}\in{\mathfrak
g}_{\alpha_i}$ and $F_{\alpha_i}\in{\mathfrak g}_{-\alpha_i}$. Also,
a basis for the Cartan subalgebra ${\mathfrak h}$ can be proposed as
$\{H_{\alpha_1}, H_{\alpha_2}, \cdots, H_{\alpha_l}\}$ such that
$\langle\alpha,H_{\alpha_i}\rangle\in\mathbb{R}$, $i=1, 2, \cdots,
l$, for all $\alpha\in\Delta$. The structure constants in form of
$\langle\beta,H_{\alpha}\rangle$ constitute the nonsingular Cartan
matrix for when $\alpha$ and $\beta$ are chosen to be the simple
roots, i.e. $a_{ij}=\langle\alpha_j,H_{\alpha_i}\rangle$, where
$a_{ii}=2$ and $a_{ij}\in\{0,-1,-2,-3\}$ for $i\neq j$. One can now
construct the full algebra starting from the Cartan-Chevalley
generators $H_{\alpha_i}$, $E_{\alpha_i}$ and $F_{\alpha_i}$, $i=1,
2, \cdots, l$, subject to the following generating relations
\begin{eqnarray}
&&\hspace{-10mm}[H_{\alpha_i},H_{\alpha_j}]=0,
\hspace{5mm}[H_{\alpha_i},E_{\alpha_j}]=a_{ij}E_{\alpha_j},
\hspace{5mm}[H_{\alpha_i},F_{\alpha_j}]=-a_{ij}F_{\alpha_j},
\hspace{5mm}[E_{\alpha_i},F_{\alpha_j}]=\delta_{i,j}H_{\alpha_i}.
\end{eqnarray}
This presentation of the Lie algebra is completed by the so-called
Serre relations
\begin{eqnarray}
&&\hspace{-10mm}\left(\mbox{ad}
E_{\alpha_i}\right)^{1-a_{ij}}E_{\alpha_j}=\sum_{k=0}^{1-a_{ij}}(-1)^k\left(\begin{array}{c}1-a_{ij}\\k\end{array}\right)
E_{\alpha_i}^{1-a_{ij}-k}E_{\alpha_j}E_{\alpha_i}^k=0,\hspace{5mm}i\neq
j,\nonumber\\\label{Serre} &&\hspace{-10mm}\left(\mbox{ad}
F_{\alpha_i}\right)^{1-a_{ij}}F_{\alpha_j}=\sum_{k=0}^{1-a_{ij}}(-1)^k\left(\begin{array}{c}1-a_{ij}\\k\end{array}\right)
F_{\alpha_i}^{1-a_{ij}-k}F_{\alpha_j}F_{\alpha_i}^k=0,\hspace{7mm}i\neq
j.
\end{eqnarray}
The structure constants in the Cartan-Chevalley basis are all
integers and they are actually the elements of the Cartan matrix.
The Lie algebra ${\mathfrak g}$ is spanned by the Cartan-Chevalley
generators and the set of multiple commutators
\begin{eqnarray}
&&\hspace{-10mm}[\cdots[[E_{\alpha_{i_{1}}},E_{\alpha_{i_{2}}}],E_{\alpha_{i_{3}}}],\cdots,E_{\alpha_{i_{k}}}],
\hspace{10mm}[\cdots[[F_{\alpha_{i_{1}}},F_{\alpha_{i_{2}}}],F_{\alpha_{i_{3}}}],\cdots,F_{\alpha_{i_{k}}}],
\end{eqnarray}
restricted by the Serre relations (\ref{Serre}). Therefore, it can
be resulted that the algebra ${\mathfrak g}$ is generated by the
Chevalley generators $E_{\alpha_i}$ and $F_{\alpha_i}$ with $i=1, 2,
\cdots, l$.

Now, let $\alpha$, $\beta\neq\pm\alpha$ and $\alpha+\beta$ be roots.
We remind the facts presented in the literature concerning the
nonzero structure constants $N_{\alpha,\beta}$  (see. e.g., Refs.
\cite{Abbaspour,Ramond}):
\begin{itemize}
\item The constants $N_{\alpha,\beta}$ are real and skew symmetric.
\item Jacobi identity:
$N_{\alpha,\beta}N_{\gamma,\eta}+N_{\beta,\gamma}N_{\alpha,\eta}+N_{\gamma,\alpha}N_{\beta,\eta}=0$,
where $\eta=-\alpha-\beta-\gamma$, and $\eta$ is not one of
$-\alpha$, $-\beta$ and $-\gamma$.
\item  If we assume the Cartan generators to be self-adjoint and $E_{\alpha}$ and $F_{\alpha}$ to be Hermitian conjugate of
each other as well as the structure constants
$\langle\beta,H_{\alpha}\rangle$ and $N_{\alpha,\beta}$ to be real
then, we will obtain
\begin{eqnarray}\label{N=-N}
N_{\alpha,\beta}=-N_{-\alpha,-\beta}. \end{eqnarray}
\item By using the properties of Killing form and the normalization condition $K(E_{\alpha},E_{-\alpha})=1$, it is shown that
$H_{\alpha}+H_{\beta}+H_{-\alpha-\beta}=0$. This relation, in turn,
leads to
\begin{eqnarray}\label{NN}
N_{\alpha,\beta}=N_{\beta,-\alpha-\beta}=N_{-\alpha-\beta,\alpha}.
\end{eqnarray}
\item Considering the $\alpha$-chain through $\beta$,
\begin{eqnarray}\label{chain}
\beta-p\alpha, \cdots, \beta-\alpha, \beta, \beta+\alpha, \cdots,
\beta+q\alpha,
\end{eqnarray}
in which, the nonnegative integer numbers $p$ and $q$ are defined by
$p-q=\langle\beta,H_{\alpha}\rangle$, we have
\begin{eqnarray}\label{NNqp}
N_{\alpha,\beta}N_{-\alpha,-\beta}=-\frac{1}{2}q(p+1)\langle\alpha,H_{\alpha}\rangle.
\end{eqnarray}
\end{itemize}

In what follows we would like to respond to the question ``whether
the structure constants in the well known commutation relations of
the exceptional Lie algebra ${\mathfrak g}_2$ saturate
simultaneously the relations (\ref{N=-N}), (\ref{NN}) and
(\ref{NNqp})". For this, in Sections 2 and 3, we consider two
different approaches to the extraction of the commutation relations
of ${\mathfrak g}_2$, where one obeys (\ref{N=-N}) and rejects
(\ref{NN}) and (\ref{NNqp}), and the other, vice versa. It is found
that the structure constants of the type of
$\langle\beta,H_{\alpha}\rangle$ are computed directly from the
Cartan matrix in the former case,  while in the latter case it does
not. In any case, the Cartan subalgebra generators are defined as
the first relation of (\ref{EEH}).

\section{The commutation relations of ${\mathfrak g}_2$ based on (\ref{N=-N}) }
The starting point is the Dynkin diagram of the rank-two exceptional
Lie algebra ${\mathfrak g}_2$:
$\bigcirc\hspace{-1.74mm}\equiv\hspace{-2mm}\equiv\hspace{-2mm}\rangle\hspace{-2mm}\equiv\hspace{-1.74mm}\bigcirc$.
The number of lines shows that the allowed angle between two roots
is $5\pi/6$. There are only two simple roots, say $\alpha_1$ and
$\alpha_2$, and we may choose them such that
$|\alpha_1|/|\alpha_2|=\sqrt{3}$. The Dynkin diagram
determines the Cartan matrix as $A=\left(\begin{array}{cc} 2 & -1 \\
-3 & 2\end{array}\right)$ (see, e.g.,
\cite{Kerf,Ramond,Georgi,Huang}). It means \setcounter{equation}{0}
\begin{eqnarray}\label{Cartannumbers}
\langle\alpha_1,H_{\alpha_1}\rangle=\langle\alpha_2,H_{\alpha_2}\rangle=2,
&\langle\alpha_1,H_{\alpha_2}\rangle=-1,&\langle\alpha_2,H_{\alpha_1}\rangle
=-3.
\end{eqnarray}
From the two last relations of (\ref{Cartannumbers}), we find that
$\alpha_1$-chain through $\alpha_2$ and the $\alpha_2$-chain through
$\alpha_1$ involve the positive roots $\alpha_2$ and
$\alpha_2+\alpha_1$ as well as $\alpha_1$, $\alpha_1+\alpha_2$,
$\alpha_1+2\alpha_2$ and $\alpha_1+3\alpha_2$, respectively. For the
last root of the second chain we have
$\langle\alpha_1+3\alpha_2,H_{\alpha_1}\rangle=-1$. This leads us to
the $\alpha_1$-chain through $\alpha_1+3\alpha_2$ with the positive
roots $\alpha_1+3\alpha_2$ and $2\alpha_1+3\alpha_2$. Then, for
${\mathfrak g}_2$, there exist six positive roots:
$\Delta_+=\{\alpha_1, \alpha_2, \alpha_1+\alpha_2,
\alpha_1+2\alpha_2, \alpha_1+3\alpha_2, 2\alpha_1+3\alpha_2\}$.
Therefore, based on the triangular decomposition, the exceptional
Lie algebra ${\mathfrak g}_2$ involves the fourteen Cartan-Weyl
generators as
\begin{eqnarray}\label{14generators}
&&H_{\alpha_1}, \ H_{\alpha_2},\ E_{\alpha_1},\ E_{\alpha_2}, \
E_{\alpha_1+\alpha_2},
\ E_{\alpha_1+2\alpha_2},\ E_{\alpha_1+3\alpha_2},\ E_{2\alpha_1+3\alpha_2},\nonumber\\
&&\hspace{22mm} F_{\alpha_1},\ F_{\alpha_2},\
F_{\alpha_1+\alpha_2},\ F_{\alpha_1+2\alpha_2}, \
F_{\alpha_1+3\alpha_2},\ F_{2\alpha_1+3\alpha_2}.
\end{eqnarray}
Now it is immediate to see that
\begin{eqnarray}\label{EFH}
&&\hspace{-8mm}[E_{\alpha_1},F_{\alpha_1}]=H_{\alpha_1},\hspace{4.2mm}
[E_{\alpha_2},F_{\alpha_2} ]=H_{\alpha_2},\\
\label{HE}
&&\hspace{-8mm}[H_{\alpha_1},E_{\alpha_1}]=2E_{\alpha_1},\hspace{2.1mm}
[H_{\alpha_1},F_{\alpha_1}]=-2F_{\alpha_1}, \hspace{2mm}
[H_{\alpha_2},E_{\alpha_2}]=2E_{\alpha_2}, \hspace{2.3mm}
[H_{\alpha_2},F_{\alpha_2}]=-2F_{\alpha_2},\\
\label{HF}
&&\hspace{-8mm}[H_{\alpha_1},E_{\alpha_2}]=-E_{\alpha_2},\
[H_{\alpha_2},E_{\alpha_1}]=-3E_{\alpha_1}, \
[H_{\alpha_1},F_{\alpha_2}]=F_{\alpha_2},\hspace{5mm}
[H_{\alpha_2},F_{\alpha_1}]=3F_{\alpha_1},\\
&&\hspace{-8mm}[H_{\alpha_1},H_{\alpha_2}]=[E_{\alpha_1},F_{\alpha_2}]=[E_{\alpha_2},F_{\alpha_1}]=0,
\end{eqnarray}
in which, nonvanishing commutation relations (\ref{HE}) and
(\ref{HF}) follow from (\ref{Cartannumbers}). With the help of the
Serre relations $(\mbox{ad}E_{\alpha_1})^{2}E_{\alpha_2}=0$,
$(\mbox{ad}E_{\alpha_2})^{4}E_{\alpha_1}=0$ and
$(\mbox{ad}E_{\alpha_1})^{2}E_{\alpha_1+3\alpha_2}=0$ one can easily
conclude that the first, the second, and the third chains of
positive roots correspond to two-, four- and two-dimensional chains
of the generators, as
\begin{eqnarray}\label{genchain1}
&&E_{\alpha_2}\stackrel{\mbox{ad}E_{\alpha_1}}{\Longrightarrow}
E_{\alpha_2+\alpha_1}
\stackrel{\mbox{ad}E_{\alpha_1}}{\Longrightarrow}0,\\\label{genchain2}
&&E_{\alpha_1}\stackrel{\mbox{ad}E_{\alpha_2}}{\Longrightarrow}
E_{\alpha_1+\alpha_2}
\stackrel{\mbox{ad}E_{\alpha_2}}{\Longrightarrow}
E_{\alpha_1+2\alpha_2}
\stackrel{\mbox{ad}E_{\alpha_2}}{\Longrightarrow}E_{\alpha_1+3\alpha_2}
\stackrel{\mbox{ad}E_{\alpha_2}}{\Longrightarrow}0
,\\\label{genchain3}
&&E_{\alpha_1+3\alpha_2}\stackrel{\mbox{ad}E_{\alpha_1}}{\Longrightarrow}
E_{2\alpha_1+3\alpha_2}
\stackrel{\mbox{ad}E_{\alpha_1}}{\Longrightarrow}0,
\end{eqnarray}
respectively. In what follows, we use the procedure proposed in Ref.
\cite{Georgi}, which is based on the use of spin representations of
the special unitary group $SU(2)$, to obtain the rest of the
commutation relations of the exceptional Lie algebra ${\mathfrak
g}_2$. According to the relations (\ref{EFH}) and (\ref{HE}), every
one of the standard triples $\{H_{\alpha_1}, E_{\alpha_1},
F_{\alpha_1}\}$ and $\{H_{\alpha_2}, E_{\alpha_2}, F_{\alpha_2}\}$
constitutes a copy of $su(2)$ commutation relations
\begin{eqnarray}\label{su(2)comm}
[J_z,J_+]=J_+,\,\,\, [J_z,J_-]=-J_-,\,\,\, [J_+,J_-]=2J_z.
\end{eqnarray}
Therefore, it is necessary to present the well known two- and
four-dimensional representations of $su(2)$ Lie algebra as below
\begin{eqnarray}
&&\hspace{-15mm}J_\pm
\left|\frac{1}{2},\mp\frac{1}{2}\right\rangle=\left|\frac{1}{2},\pm\frac{1}{2}\right\rangle,\,\,\,\,\,\,\,\,\,\,\,
J_\pm\left|\frac{1}{2},\pm\frac{1}{2}\right\rangle=0,\\
&&\hspace{-15mm}J_\pm \left|
\frac{3}{2},\mp\frac{3}{2}\right\rangle=\sqrt{3}\left|\frac{3}{2},\mp\frac{1}{2}\right\rangle,
\ J_\pm
\left|\frac{3}{2},\mp\frac{1}{2}\right\rangle=2\left|\frac{3}{2},\pm\frac{1}{2}\right\rangle,
\nonumber\\ &&\hspace{-15mm}
J_\pm\left|\frac{3}{2},\pm\frac{1}{2}\right\rangle=\sqrt{3}\left|\frac{3}{2},
\pm\frac{3}{2}\right\rangle,\ J_\pm\left|\frac{3}{2},\pm
\frac{3}{2}\right\rangle=0.
\end{eqnarray}
Comparing the relations (\ref{EFH}) and (\ref{HE}) with
(\ref{su(2)comm}), it is found that the two-, four- and
two-dimensional chains (\ref{genchain1}), (\ref{genchain2}) and
(\ref{genchain3}) for the generators can be corresponded to the
following definitions
\begin{eqnarray}
&&\hspace{-5mm}E_{\alpha_2}\equiv\left|\frac{1}{2},\frac{-1}{2};1\right\rangle,\
E_{\alpha_1+\alpha_2}\equiv\left|\frac{1}{2},\frac{1}{2};1\right\rangle,
\ \mbox{ad}H_{\alpha_1}\equiv2J_z, \
\mbox{ad}E_{\alpha_1}\equiv-\sqrt{\frac{3}{2}}J_+,\
\mbox{ad}F_{\alpha_1}\equiv-\sqrt{\frac{2}{3}}J_-,\nonumber\\\label{rep1}\\
&&\hspace{-5mm}E_{\alpha_1}\equiv\left|\frac{3}{2},\frac{-3}{2};2\right\rangle,
\
E_{\alpha_1+\alpha_2}\equiv\left|\frac{3}{2},\frac{-1}{2};2\right\rangle,
\
E_{\alpha_1+2\alpha_2}\equiv\left|\frac{3}{2},\frac{1}{2};2\right\rangle,
\
E_{\alpha_1+3\alpha_2}\equiv\left|\frac{3}{2},\frac{3}{2};2\right\rangle,
\nonumber\\
&&\hspace{-5mm}\mbox{ad}H_{\alpha_2}\equiv2J_z,\
\mbox{ad}E_{\alpha_2}\equiv\frac{J_+}{\sqrt{2}},\
\mbox{ad}F_{\alpha_2}\equiv\sqrt{2} J_-, \label{rep2}\\
&&\hspace{-5mm}E_{\alpha_1+3\alpha_2}\equiv\left|\frac{1}{2},\frac{-1}{2};3\right\rangle,
\
E_{2\alpha_1+3\alpha_2}\equiv\left|\frac{1}{2},\frac{1}{2};3\right\rangle,
\ \mbox{ad}H_{\alpha_1}\equiv2J_z,\
\mbox{ad}E_{\alpha_1}\equiv\sqrt{\frac{3}{2}}J_+,\
\mbox{ad}F_{\alpha_1}\equiv\sqrt{\frac{2}{3}}J_-, \nonumber\\
\label{rep3}
\end{eqnarray}
respectively. Here, the labels $1$, $2$ and $3$ on the right end of
the kets denote the chain issue. From now on, for simplicity, the
skew symmetric property  is implicitly used to determine the
structure constants. From (\ref{rep2}) and (\ref{rep3}) we
immediately get
\begin{eqnarray}\label{Ea1a2}
&&\hspace{-5mm}\mbox{ad}E_{\alpha_2}E_{\alpha_1}\equiv\frac{J_+}{\sqrt{2}}\left|\frac{3}{2},
\frac{-3}{2},2\right\rangle=\sqrt{\frac{3}{2}}\left|\frac{3}{2},\frac{-1}{2},2\right\rangle
\equiv\sqrt{\frac{3}{2}}E_{\alpha_1+\alpha_2} \Rightarrow
[E_{\alpha_2},E_{\alpha_1}]=\sqrt{\frac{3}{2}}E_{\alpha_1+\alpha_2},\\
&&\hspace{-5mm}\mbox{ad}E_{\alpha_2}E_{\alpha_1+\alpha_2}\equiv\frac{J_+}{\sqrt{2}}\left|\frac{3}{2},\frac{-1}{2},2\right\rangle
=\sqrt{2}\left|\frac{3}{2},\frac{1}{2},2\right\rangle\equiv\sqrt{2}E_{\alpha_1+2\alpha_2}
\Rightarrow[E_{\alpha_2},E_{\alpha_1+\alpha_2}]=\sqrt{2}
E_{\alpha_1+2\alpha_2},\nonumber\\ \\
&&\hspace{-5mm}\mbox{ad}E_{\alpha_2}E_{\alpha_1+2\alpha_2}\equiv\frac{J_+}{\sqrt{2}}\left|\frac{3}{2},\frac{1}{2},2\right\rangle
=\sqrt{\frac{3}{2}}\left|\frac{3}{2},\frac{3}{2},2\right\rangle\equiv\sqrt{\frac{3}{2}}E_{\alpha_1+3\alpha_2}
\Rightarrow[E_{\alpha_2},E_{\alpha_1+2\alpha_2}]=\sqrt{\frac{3}{2}}
E_{\alpha_1+3\alpha_2},\nonumber\\ \\
&&\hspace{-5mm}\mbox{ad}E_{\alpha_1}
E_{\alpha_1+3\alpha_2}\equiv\sqrt{\frac{3}{2}}
J_+\left|\frac{1}{2},\frac{-1}{2};3\right\rangle
=\sqrt{\frac{3}{2}}\left|\frac{1}{2},\frac{1}{2};3\right\rangle\equiv\sqrt{\frac{3}{2}}E_{2\alpha_1+3\alpha_2}\nonumber\\
&&\hspace{85mm}
\Rightarrow[E_{\alpha_1},E_{\alpha_1+3\alpha_2}]=\sqrt{\frac{3}{2}}
E_{2\alpha_1+3\alpha_2}.
\end{eqnarray}
Note that the result (\ref{Ea1a2}) can be also obtained by using
(\ref{rep1}). The condition (\ref{N=-N}) using the latter four
relations will enable us to determine the following four commutation
relations
\begin{eqnarray}
\begin{array}{lll}
&[F_{\alpha_2 },F_{\alpha_1}]=-\sqrt{\frac{3}{2}}
F_{\alpha_1+\alpha_2},
& [F_{\alpha_2},F_{\alpha_1+\alpha_2}]=-\sqrt{2} F_{\alpha_1+2\alpha_2},\\
&[F_{\alpha_2 },F_{\alpha_1+2\alpha_2}]=-\sqrt{\frac{3}{2}}
F_{\alpha_1+3\alpha_2}, &
[F_{\alpha_1},F_{\alpha_1+3\alpha_2}]=-\sqrt{\frac{3}{2}}
F_{2\alpha_1+3\alpha_2}.
\end{array}
\end{eqnarray}
The next commutation relations follow easily from the fact that the
number of the positive roots of ${\mathfrak g}_2$ is exactly six:
\begin{eqnarray}
\begin{array}{ll}
&[E_{\alpha_1},E_{\alpha_1+\alpha_2}]=[E_{\alpha_1},E_{\alpha_1+2\alpha_2}]=
[E_{\alpha_1},E_{2\alpha_1+3\alpha_2}]=[E_{\alpha_2},E_{\alpha_1+3\alpha_2}]=0,\\
&[E_{\alpha_2},E_{2\alpha_1+3\alpha_2}]=[E_{\alpha_1+\alpha_2},E_{\alpha_1+3\alpha_2}]=
[E_{\alpha_1+\alpha_2},E_{2\alpha_1+3\alpha_2}]=0,\\
&[E_{\alpha_1+2\alpha_2},E_{2\alpha_1+3\alpha_2}]=[E_{\alpha_1+3\alpha_2},E_{2\alpha_1+3\alpha_2}]=
[E_{\alpha_1+2\alpha_2},E_{\alpha_1+3\alpha_2}]=0.
\end{array}
\end{eqnarray}
Similar relations are also held for the $F$-generators corresponding
to the negative roots. The above relations, together with the Jacobi
identity, give rise to
\begin{eqnarray}
\begin{array}{lll}
[E_{\alpha_1+\alpha_2},E_{\alpha_1+2\alpha_2}]&=&-\sqrt{\frac{2}{3}}[[E_{\alpha_1},E_{\alpha_2}],E_{\alpha_1+2\alpha_2}]\\
&=&\sqrt{\frac{2}{3}}[[E_{\alpha_1+2\alpha_2},E_{\alpha_1}],E_{_2
}]+
\sqrt{\frac{2}{3}}[[E_{\alpha_2},E_{\alpha_1+2\alpha_2}],E_{\alpha_1}]=-\sqrt{\frac{3}{2}}
E_{2\alpha_1+3\alpha_2}.
\end{array}
\end{eqnarray}
In the same way, from the Jacobi identity and earlier commutation
relations we find
\begin{eqnarray}
\begin{array}{lllll}
&[H_{\alpha_1},E_{\alpha_1+\alpha_2}]=E_{\alpha_1+\alpha_2},&
[H_{\alpha_1},E_{\alpha_1+2\alpha_2}]=0,&
[H_{\alpha_1},E_{\alpha_1+3\alpha_2}]=-E_{\alpha_1+3\alpha_2},\\
&[H_{\alpha_1},E_{2\alpha_1+3\alpha_2}]=E_{2\alpha_1+3\alpha_2},&
[H_{\alpha_2},E_{\alpha_1+\alpha_2}]=-E_{\alpha_1+\alpha_2},&
[H_{\alpha_2},E_{\alpha_1+2\alpha_2}]=E_{\alpha_1+2\alpha_2},\\
&[H_{\alpha_2},E_{\alpha_1+3\alpha_2}]=3E_{\alpha_1+3\alpha_2},&
[H_{\alpha_2},E_{2\alpha_1+3\alpha_2}]=0,&
\end{array}
\end{eqnarray}
\begin{eqnarray}\label{HEFHH}
\begin{array}{ll}
&H_{\alpha_1+\alpha_2}\equiv[E_{\alpha_1+\alpha_2},F_{\alpha_1+\alpha_2}]=2H_{\alpha_1}+\frac{2}{3}H_{\alpha_2},
\\ &H_{\alpha_1+2\alpha_2}\equiv[E_{\alpha_1+2\alpha_2},F_{\alpha_1+2\alpha_2}]=4H_{\alpha_1}+\frac{8}{3} H_{\alpha_2},\\
&H_{\alpha_1+3\alpha_2}\equiv[E_{\alpha_1+3\alpha_2},F_{\alpha_1+3\alpha_2}]=8H_{\alpha_1}+8H_{\alpha_2},\\
&H_{2\alpha_1+3\alpha_2}\equiv[E_{2\alpha_1+3\alpha_2},F_{2\alpha_1+3\alpha_2}]=\frac{16}{3}(2H_{\alpha_1}+H_{\alpha_2}),
\end{array}
\end{eqnarray}
\begin{eqnarray}
\begin{array}{lll}
&[H_{\alpha_1+\alpha_2},E_{\alpha_1+\alpha_2}]=\frac{4}{3}E_{\alpha_1+\alpha_2},&
[H_{\alpha_1+\alpha_2},E_{\alpha_1+2\alpha_2}]=\frac{2}{3}E_{\alpha_1+2\alpha_2},\\
&[H_{\alpha_1+\alpha_2},E_{\alpha_1+3\alpha_2}]=0,&
[H_{\alpha_1+\alpha_2},E_{2\alpha_1+3\alpha_2}]=2E_{2\alpha_1+3\alpha_2},\\
&[H_{\alpha_1+2\alpha_2},E_{\alpha_1+\alpha_2}]=\frac{4}{3}E_{\alpha_1+\alpha_2},&
[H_{\alpha_1+2\alpha_2},E_{\alpha_1+2\alpha_2}]=\frac{8}{3}E_{\alpha_1+2\alpha_2},\\
&[H_{\alpha_1+2\alpha_2},E_{_1+3\alpha_2}]=4E_{\alpha_1+3\alpha_2},&
[H_{\alpha_1+2\alpha_2},E_{2\alpha_1+3\alpha_2}]=4E_{2\alpha_1+3\alpha_2},\\
&[H_{\alpha_1+3\alpha_2},E_{\alpha_1+\alpha_2}]=0,&
[H_{\alpha_1+3\alpha_2},E_{\alpha_1+2\alpha_2}]=8E_{\alpha_1+2\alpha_2},\\
&[H_{\alpha_1+3\alpha_2},E_{\alpha_1+3\alpha_2}]=16E_{\alpha_1+3\alpha_2},&
[H_{\alpha_1+3\alpha_2},E_{2\alpha_1+3\alpha_2}]=8E_{2\alpha_1+3\alpha_2},\\
&[H_{2\alpha_1+3\alpha_2},E_{\alpha_1+\alpha_2}]=\frac{16}{3}E_{\alpha_1+\alpha_2},&
[H_{2\alpha_1+3\alpha_2},E_{\alpha_1+2\alpha_2}]=\frac{16}{3}E_{\alpha_1+2\alpha_2},\\
&[H_{2\alpha_1+3\alpha_2},E_{\alpha_1+3\alpha_2}]=\frac{16}{3}E_{\alpha_1+3\alpha_2},&
[H_{2\alpha_1+3\alpha_2},E_{2\alpha_1+3\alpha_2}]=\frac{32}{3}E_{2\alpha_1+3\alpha_2},
\end{array}
\end{eqnarray}
\begin{eqnarray}
\begin{array}{llll}
&[F_{\alpha_1},E_{\alpha_1+\alpha_2}]=-\sqrt{\frac{2}{3}}E_{\alpha_2},
&[F_{\alpha_1},E_{\alpha_1+2\alpha_2}]=0, &
[F_{\alpha_1},E_{\alpha_1+3\alpha_2}]=0,
\\
&[F_{\alpha_1},E_{2\alpha_1+3\alpha_2}]=\sqrt{\frac{2}{3}}E_{\alpha_1+3\alpha_2},
& [F_{\alpha_2},E_{\alpha_1+\alpha_2}]=\sqrt{6}E_{\alpha_1}, &
[F_{\alpha_2},E_{_1+2\alpha_2}]=2\sqrt{2}E_{\alpha_1+\alpha_2},
\\
&[F_{\alpha_2},E_{\alpha_1+3\alpha_2}]=\sqrt{6}E_{_1+2\alpha_2},
&[F_{\alpha_2},E_{2\alpha_1+3\alpha_2}]=0. &
\end{array}
\end{eqnarray}
Now, with the help of (\ref{HEFHH}) we can directly calculate the
following commutation relations without the use of the Jacobi
identity
\begin{eqnarray}
\begin{array}{llll}
&[H_{\alpha_1+\alpha_2},E_{\alpha_1}]=2E_{\alpha_1},&
[H_{\alpha_1+\alpha_2},E_{\alpha_2}]=\frac{-2}{3}E_{\alpha_2},& [H_{\alpha_1+2\alpha_2},E_{\alpha_1}]=0,\\
&[H_{\alpha_1+2\alpha_2},E_{\alpha_2}]=\frac{4}{3}E_{\alpha_2},&
[H_{\alpha_1+3\alpha_2},E_{\alpha_1}]=-8E_{\alpha_1},&
[H_{\alpha_1+3\alpha_2},E_{\alpha_2}]=8E_{\alpha_2}, \\
&[H_{2\alpha_1+3\alpha_2},E_{\alpha_1}]=\frac{16}{3}E_{\alpha_1},&
[H_{2\alpha_1+3\alpha_2},E_{\alpha_2}]=0.&
\end{array}
\end{eqnarray}
The condition (\ref{N=-N}) together with the Jacobi identity, leads
to
\begin{eqnarray}
\begin{array}{lll}
&[E_{\alpha_1+\alpha_2},F_{\alpha_1+2\alpha_2}]=\frac{4\sqrt{2}}{3}
F_{\alpha_2}, & [E_{\alpha_1+\alpha_2},F_{\alpha_1+3\alpha_2}]=0, \\
&[E_{\alpha_1+\alpha_2},F_{2\alpha_1+3\alpha_2}]=2\sqrt{\frac{2}{3}}F_{\alpha_1+2\alpha_2},
&[E_{\alpha_1+2\alpha_2},F_{\alpha_1+\alpha_2}]=\frac{4\sqrt{2}}{3}E_{\alpha_2},\\
&[E_{\alpha_1+2\alpha_2},F_{\alpha_1+3\alpha_2}]=4\sqrt{\frac{2}{3}}F_{\alpha_2},
&[E_{\alpha_1+2\alpha_2},F_{2\alpha_1+3\alpha_2}]=-4\sqrt{\frac{2}{3}}F_{\alpha_1+\alpha_2},\\
&[E_{\alpha_1+3\alpha_2},F_{\alpha_1+\alpha_2}]=0,&
[E_{\alpha_1+3\alpha_2},F_{\alpha_1+2\alpha_2}]=4\sqrt{\frac{2}{3}}E_{\alpha_2},
\\
&[E_{\alpha_1+3\alpha_2},F_{2\alpha_1+3\alpha_2}]=8\sqrt{\frac{2}{3}}F_{\alpha_1},
&[E_{2\alpha_1+3\alpha_2},F_{\alpha_1+\alpha_2}]=2\sqrt{\frac{2}{3}}E_{\alpha_1+2\alpha_2},\\
&[E_{2\alpha_1+3\alpha_2},F_{\alpha_1+2\alpha_2}]=-4\sqrt{\frac{2}{3}}E_{\alpha_1+\alpha_2},
&[E_{2\alpha_1+3\alpha_2},F_{\alpha_1+3\alpha_2}]=8\sqrt{\frac{2}{3}}E_{\alpha_1},\\
&[E_{\alpha_1},F_{\alpha_1+\alpha_2}]=\sqrt{\frac{2}{3}}F_{\alpha_2},
&
[E_{\alpha_1},F_{\alpha_1+2\alpha_2}]=0, \\
&[E_{\alpha_1},F_{\alpha_1+3\alpha_2}]=0,&
[E_{\alpha_1},F_{2\alpha_1+3\alpha_2}]=-\sqrt{\frac{2}{3}}F_{\alpha_1+3\alpha_2},\\
&[E_{\alpha_2},F_{\alpha_1+\alpha_2}]=-\sqrt{6}F_{\alpha_1}, &
[E_{\alpha_2},F_{\alpha_1+2\alpha_2}]=-2\sqrt{2}F_{\alpha_1+\alpha_2},\\
&[E_{\alpha_2},F_{\alpha_1+3\alpha_2}]=-\sqrt{6}F_{\alpha_1+2\alpha_2},&
[E_{\alpha_2},F_{2\alpha_1+3\alpha_2}]=0, \\
&[F_{\alpha_1},F_{\alpha_2}]=\sqrt{\frac{3}{2}}F_{\alpha_1+\alpha_2},&
[F_{\alpha_1},F_{\alpha_1+3\alpha_2}]=-\sqrt{\frac{3}{2}}F_{2\alpha_1+3\alpha_2},\\
&[F_{\alpha_2},F_{\alpha_1+\alpha_2}]=-\sqrt{2}F_{\alpha_1+2\alpha_2},&
[F_{\alpha_2},F_{\alpha_1+2\alpha_2}]=-\sqrt{\frac{3}{2}}F_{\alpha_1+3\alpha_2},\\
&[F_{\alpha_1+\alpha_2},F_{\alpha_1+2\alpha_2}]=\sqrt{\frac{3}{2}}F_{2\alpha_1+3\alpha_2}.&
\end{array}
\end{eqnarray}
Now, we can calculate the commutator of elements of the Cartan
subalgebra with the $F$-generators corresponding to a positive and
non-simple root. For example,
\begin{eqnarray}
[H_{\alpha_2},F_{\alpha_1+\alpha_2}]&=&\sqrt{\frac{2}{3}}[H_{\alpha_2},[F_{\alpha_1},F_{\alpha_2}]]
=-\sqrt{\frac{2}{3}}[F_{\alpha_2},[H_{\alpha_2},F_{\alpha_1}]]-\sqrt{\frac{2}{3}}[F_{\alpha_1},[F_{\alpha_2},H_{\alpha_2}]]
\nonumber\\
&=&-\sqrt{\frac{2}{3}}[F_{\alpha_2},F_{\alpha_1}]=F_{\alpha_1+\alpha_2}.
\end{eqnarray}
By rescaling the generators corresponding to the positive non-simple
roots as
\begin{eqnarray}
\begin{array}{llllll}
& X_{\alpha_1}\equiv E_{\alpha_1}, & X_{\alpha_2}\equiv
E_{\alpha_2}, &  X_{\alpha_1+\alpha_2}\equiv
-\sqrt{\frac{3}{2}}E_{\alpha_1+\alpha_2}, \\ &
X_{\alpha_1+2\alpha_2}\equiv-\frac{\sqrt{3}}{2}E_{\alpha_1+2\alpha_2},
&X_{\alpha_1+3\alpha_2}\equiv-\frac{\sqrt{2}}{4}
E_{\alpha_1+3\alpha_2}, &
X_{2\alpha_1+3\alpha_2}\equiv-\frac{\sqrt{3}}{4}
E_{2\alpha_1+3\alpha_2}, & &
\end{array}
\end{eqnarray}
and likewise for the negative roots denoted by $Y$, we obtain the
familiar commutation relations of the exceptional Lie algebra
${\mathfrak g}_2$, as shown in  Table 1 (see e.g., Refs.
\cite{Saenkarun1,Huang,Wildberger,Saenkarun2}).

\section{The commutation relations of ${\mathfrak g}_2$ based on (\ref{NN}) and
(\ref{NNqp})} In this section we will denote the Cartan-Chevalley
generators by $H^{\prime}_{\alpha_i}$, $X^{\prime}_{\alpha_i}$ and
$Y^{\prime}_{\alpha_i}$ which correspond respectively to the zero
and the simple (positive and negative) roots. The conditions
(\ref{NN}) for the root system of the exceptional Lie algebra
${\mathfrak g}_2$ can be explicitly described as follows
\setcounter{equation}{0}
\begin{eqnarray}\label{NNqpg2}
\begin{array}{lll}
&N_{\alpha_1,\alpha_2}=N_{\alpha_2,-\alpha_1-\alpha_2}=N_{-\alpha_1-\alpha_2,\alpha_1},\\
&N_{-\alpha_1,-\alpha_2}=N_{-\alpha_2,\alpha_1+\alpha_2}=N_{\alpha_1+\alpha_2,-\alpha_1},
\\
&N_{\alpha_1,\alpha_1+3\alpha_2}=N_{\alpha_1+3\alpha_2,-2\alpha_1-3\alpha_2}=
N_{-2\alpha_1-3\alpha_2,\alpha_1},\\
&N_{-\alpha_1,-\alpha_1-3\alpha_2}=N_{-\alpha_1-3\alpha_2,2\alpha_1+3\alpha_2}=N_{2\alpha_1+3\alpha_2,-\alpha_1},
\\
&N_{\alpha_2,\alpha_1+\alpha_2}=N_{\alpha_1+\alpha_2,-\alpha_1-2\alpha_2}=N_{-\alpha_1-2\alpha_2,\alpha_2},\\
&N_{-\alpha_2,-\alpha_1-\alpha_2}=N_{-\alpha_1-\alpha_2,\alpha_1+2\alpha_2}=N_{\alpha_1+2\alpha_2,-\alpha_2},
\\
&N_{\alpha_2,\alpha_1+2\alpha_2}=N_{\alpha_1+2\alpha_2,-\alpha_1-3\alpha_2}=N_{-\alpha_1-3\alpha_2,\alpha_2},\\
&N_{-\alpha_2,-\alpha_1-2\alpha_2}=N_{-\alpha_1-2\alpha_2,\alpha_1+3\alpha_2}=N_{\alpha_1+3\alpha_2,-\alpha_2},
\\
&N_{\alpha_1+\alpha_2,\alpha_1+2\alpha_2}=N_{\alpha_1+2\alpha_2,-2\alpha_1-3\alpha_2}=N_{-2\alpha_1-3\alpha_2,\alpha_1+\alpha_2},\\
&N_{-\alpha_1-\alpha_2,-\alpha_1-2\alpha_2}=N_{-\alpha_1-2\alpha_2,2\alpha_1+3\alpha_2}=N_{2\alpha_1+3\alpha_2,-\alpha_1-\alpha_2}.
\end{array}
\end{eqnarray}
One can investigate and see that there exist fifteen various root
chains of the type (\ref{chain}),  each of them corresponds to one
relation in the form of (\ref{NNqp}):
\begin{eqnarray}\label{chaing2}
\begin{array}{l}
\left\{
\begin{array}{l}
\beta=\alpha_2,\,\alpha=\alpha_1,\,p=0,\,q=1\hspace{3mm}\&\hspace{3mm}\beta=\alpha_1,\,\alpha=\alpha_2,\,p=0,\,q=3,\\
N_{\alpha_1,\alpha_2}N_{-\alpha_1,-\alpha_2}=\frac{-1}{2}\langle
\alpha_1,H'_{\alpha_1}\rangle= \frac{-3}{2} \langle
\alpha_2,H'_{\alpha_2}\rangle,
\end{array} \right.\\ \\
\left\{
\begin{array}{l}
\beta=\alpha_1+3\alpha_2,\,\alpha=\alpha_1,\,p=0,\,q=1\hspace{3mm}\&\hspace{3mm}\beta=\alpha_1,\,\alpha=\alpha_1+3\alpha_2,\,p=0,\,q=1,\\
N_{\alpha_1,\alpha_1+3\alpha_2} N_{-\alpha_1,-\alpha_1-3\alpha_2}
=\frac{-1}{2} \langle \alpha_1,H'_{\alpha_1}\rangle =\frac{-1}{2}
\langle \alpha_1+3\alpha_2,H'_{\alpha_1+3\alpha_2}\rangle,
\end{array} \right.\\ \\
\left\{
\begin{array}{l}
\beta=-\alpha_1-\alpha_2,\,\alpha=\alpha_1,\,p=0,\,q=1\hspace{3mm}\&\hspace{3mm}\beta=\alpha_1,\,\alpha=-\alpha_1-\alpha_2,\,p=0,\,q=3,\\
N_{\alpha_1,-\alpha_1-\alpha_2} N_{-\alpha_1,\alpha_1+\alpha_2}
=\frac{-1}{2} \langle \alpha_1,H'_{\alpha_1}\rangle =\frac{-3}{2}
\langle \alpha_1+\alpha_2,H'_{\alpha_1+\alpha_2}\rangle,
\end{array} \right.\\ \\
\left\{
\begin{array}{l}
\beta=-2\alpha_1-3\alpha_2,\,\alpha=\alpha_1,\,p=0,\,q=1\hspace{3mm}\&\hspace{3mm}\beta=\alpha_1,\,\alpha=-2\alpha_1-3\alpha_2,\,p=0,\,q=1,\\
N_{\alpha_1,-2\alpha_1-3\alpha_2}N_{-\alpha_1,2\alpha_1+3\alpha_2}=\frac{-1}{2}\langle
\alpha_1,H'_{\alpha_1}\rangle= \frac{-1}{2} \langle
2\alpha_1+3\alpha_2,H'_{2\alpha_1+3\alpha_2}\rangle,
\end{array} \right.\\ \\
\left\{
\begin{array}{l}
\beta=\alpha_1+\alpha_2,\,\alpha=\alpha_2,\,p=1,\,q=2\hspace{3mm}\&\hspace{3mm}\beta=\alpha_2,\,\alpha=\alpha_1+\alpha_2,\,p=1,\,q=2,\\
N_{\alpha_2,\alpha_1+\alpha_2}N_{-\alpha_2,-\alpha_1-\alpha_2}=-2\langle
\alpha_2,H'_{\alpha_2}\rangle =-2\langle
\alpha_1+\alpha_2,H'_{\alpha_1+\alpha_2}\rangle,
\end{array} \right.\\ \\
\left\{
\begin{array}{l}
\beta=\alpha_1+2\alpha_2,\,\alpha=\alpha_2,\,p=2,\,q=1\hspace{3mm}\&\hspace{3mm}\beta=\alpha_2,\,\alpha=\alpha_1+2\alpha_2,\,p=2,\,q=1,\\
N_{\alpha_2,\alpha_1+2\alpha_2}N_{-\alpha_2,-\alpha_1-2\alpha_2}=\frac{-3}{2}\langle
\alpha_2,H'_{\alpha_2}\rangle= \frac{-3}{2} \langle
\alpha_1+2\alpha_2,H'_{\alpha_1+2\alpha_2}\rangle,
\end{array} \right.\\ \\
\left\{
\begin{array}{l}
\beta=-\alpha_1-3\alpha_2,\,\alpha=\alpha_2,\,p=0,\,q=3\hspace{3mm}\&\hspace{3mm}\beta=\alpha_2,\,\alpha=-\alpha_1-3\alpha_2,\,p=0,\,q=1,\\
N_{\alpha_2,-\alpha_1-3\alpha_2}N_{-\alpha_2,\alpha_1+3\alpha_2}=\frac{-3}{2}\langle
\alpha_2,H'_{\alpha_2}\rangle= \frac{-1}{2} \langle
\alpha_1+3\alpha_2,H'_{\alpha_1+3\alpha_2}\rangle,
\end{array} \right.\\ \\
\left\{
\begin{array}{l}
\beta=\alpha_1+2\alpha_2,\,\alpha=\alpha_1+\alpha_2,\,p=2,\,q=1\hspace{3mm}\&\hspace{3mm}\beta=\alpha_1+\alpha_2,\,\alpha=\alpha_1+2\alpha_2,\,p=2,\,q=1,\\
N_{\alpha_1+\alpha_2,\alpha_1+2\alpha_2}N_{-\alpha_1-\alpha_2,-\alpha_1-2\alpha_2}
=\frac{-3}{2}\langle
\alpha_1+\alpha_2,H'_{\alpha_1+\alpha_2}\rangle=\frac{-3}{2} \langle
\alpha_1+2\alpha_2,H'_{\alpha_1+2\alpha_2}\rangle,
\end{array} \right.\\ \\
\left\{
\begin{array}{l}
\beta=-2\alpha_1-3\alpha_2,\,\alpha=\alpha_1+\alpha_2,\,p=0,\,q=3\hspace{3mm}\&\hspace{3mm}\beta=\alpha_1+\alpha_2,\,\alpha=-2\alpha_1-3\alpha_2,\,p=0,\,q=1,\\
N_{\alpha_1+\alpha_2,-2\alpha_1-3\alpha_2}N_{-\alpha_1-\alpha_2,2\alpha_1+3\alpha_2}
=\frac{-3}{2}\langle
\alpha_1+\alpha_2,H'_{\alpha_1+\alpha_2}\rangle=\frac{-1}{2} \langle
2\alpha_1+3\alpha_2,H'_{2\alpha_1+3\alpha_2}\rangle,
\end{array} \right.\\ \\
\left\{
\begin{array}{l}
\beta=-\alpha_1-3\alpha_2,\,\alpha=\alpha_1+2\alpha_2,\,p=0,\,q=3\hspace{3mm}\&\hspace{3mm}\beta=\alpha_1+2\alpha_2,\,\alpha=-\alpha_1-3\alpha_2,\,p=0,\,q=1,\\
N_{\alpha_1+2\alpha_2,-\alpha_1-3\alpha_2}N_{-\alpha_1-2\alpha_2,\alpha_1+3\alpha_2}
=\frac{-3}{2}\langle
\alpha_1+2\alpha_2,H'_{\alpha_1+2\alpha_2}\rangle=\frac{-1}{2}
\langle \alpha_1+3\alpha_2,H'_{\alpha_1+3\alpha_2}\rangle,
\end{array} \right.\\ \\
\left\{
\begin{array}{l}
\beta=-2\alpha_1-3\alpha_2,\,\alpha=\alpha_1+3\alpha_2,\,p=0,\,q=1\hspace{3mm}\&\hspace{3mm}\beta=\alpha_1+3\alpha_2,\,\alpha=-2\alpha_1-3\alpha_2,\,p=0,\,q=1,\\
N_{\alpha_1+3\alpha_2,-2\alpha_1-3\alpha_2}N_{-\alpha_1-3\alpha_2,2\alpha_1+3\alpha_2}
=\frac{-1}{2}\langle
\alpha_1+3\alpha_2,H'_{\alpha_1+3\alpha_2}\rangle=\frac{-1}{2}
\langle 2\alpha_1+3\alpha_2,H'_{2\alpha_1+3\alpha_2}\rangle,
\end{array} \right.\\ \\
\left\{
\begin{array}{l}
\beta=\alpha_1+3\alpha_2,\,\alpha=-\alpha_1-2\alpha_2,\,p=0,\,q=3\hspace{3mm}\&\hspace{3mm}\beta=-\alpha_1-2\alpha_2,\,\alpha=\alpha_1+3\alpha_2,\,p=0,\,q=1,\\
N_{-\alpha_1-2\alpha_2,\alpha_1+3\alpha_2}N_{\alpha_1+2\alpha_2,-\alpha_1-3\alpha_2}
=\frac{-3}{2}\langle
\alpha_1+2\alpha_2,H'_{\alpha_1+2\alpha_2}\rangle=\frac{-1}{2}
\langle \alpha_1+3\alpha_2,H'_{\alpha_1+3\alpha_2}\rangle,
\end{array} \right.\\ \\
\left\{
\begin{array}{l}
\beta=2\alpha_1+3\alpha_2,\,\alpha=-\alpha_1-3\alpha_2,\,p=0,\,q=1\hspace{3mm}\&\hspace{3mm}\beta=-\alpha_1-3\alpha_2,\,\alpha=2\alpha_1+3\alpha_2,\,p=0,\,q=1,\\
N_{-\alpha_1-3\alpha_2,2\alpha_1+3\alpha_2}N_{\alpha_1+3\alpha_2,-2\alpha_1-3\alpha_2}
=\frac{-1}{2}\langle
\alpha_1+3\alpha_2,H'_{\alpha_1+3\alpha_2}\rangle=\frac{-1}{2}
\langle 2\alpha_1+3\alpha_2,H'_{2\alpha_1+3\alpha_2}\rangle,
\end{array} \right.\\ \\
\left\{
\begin{array}{l}
\beta=-\alpha_1-\alpha_2,\,\alpha=\alpha_2,\,p=2,\,q=1\hspace{3mm}\&\hspace{3mm}\beta=\alpha_2,\,\alpha=-\alpha_1-\alpha_2,\,p=2,\,q=1,\\
N_{\alpha_2,-\alpha_1-\alpha_2}N_{-\alpha_2,\alpha_1+\alpha_2}
=\frac{-3}{2}\langle \alpha_2,H'_{\alpha_2}\rangle=\frac{-3}{2}
\langle \alpha_1+\alpha_2,H'_{\alpha_1+\alpha_2}\rangle,
\end{array} \right.\\ \\
\left\{
\begin{array}{l}
\beta=-\alpha_1-2\alpha_2,\,\alpha=\alpha_2,\,p=1,\,q=2\hspace{3mm}\&\hspace{3mm}\beta=\alpha_2,\,\alpha=-\alpha_1-2\alpha_2,\,p=1,\,q=2,\\
N_{\alpha_2,-\alpha_1-2\alpha_2}N_{-\alpha_2,\alpha_1+2\alpha_2}
=-2\langle \alpha_2,H'_{\alpha_2}\rangle=-2\langle
\alpha_1+2\alpha_2,H'_{\alpha_1+2\alpha_2}\rangle,
\end{array} \right.\\ \\
\end{array}
\end{eqnarray}
Considering the relations (\ref{NNqpg2}) and (\ref{chaing2}),
without loss of generality, we can fix the structure constants in
the form of $\langle\beta,H_{\alpha}^{\prime}\rangle$ as follows
\begin{eqnarray}\label{one}
%&&\hspace{-20mm}
\langle2\alpha_1+3\alpha_2,H'_{2\alpha_1+3\alpha_2}\rangle&=&\langle
\alpha_1+3\alpha_2,H'_{\alpha_1+3\alpha_2}\rangle =\langle
\alpha_1,H'_{\alpha_1}\rangle=3\langle
\alpha_1+\alpha_2,H'_{\alpha_1+\alpha_2}\rangle
\nonumber\\
&=&3\langle \alpha_1+2\alpha_2,
H'_{\alpha_1+2\alpha_2}\rangle=3\langle
\alpha_2,H'_{\alpha_2}\rangle=1.
\end{eqnarray}
Substituting the results (\ref{one}) into equations (\ref{NNqpg2})
and (\ref{chaing2}), we obtain the following values for the
structure constants in the form of $N_{\alpha,\beta}$,
\begin{eqnarray}\label{Ng21}
N_{\alpha_1,\alpha_2}&=&N_{\alpha_1,\alpha_1+3\alpha_2}=-N_{\alpha_1,-\alpha_1-\alpha_2}=-N_{\alpha_1,-2\alpha_1-3\alpha_2}
=N_{\alpha_2,\alpha_1+\alpha_2}=N_{\alpha_2,\alpha_1+2\alpha_2}
\nonumber\\
&=&N_{\alpha_2,-\alpha_1-\alpha_2}=-N_{\alpha_2,-\alpha_1-2\alpha_2}
=-N_{\alpha_2,-\alpha_1-3\alpha_2}=N_{\alpha_1+\alpha_2,\alpha_1+2\alpha_2}
\nonumber\\
&=&N_{\alpha_1+\alpha_2,-\alpha_1-2\alpha_2}
=-N_{\alpha_1+\alpha_2,-2\alpha_1-3\alpha_2}=N_{\alpha_1+2\alpha_2,-\alpha_1-3\alpha_2}
=N_{\alpha_1+2\alpha_2,-2\alpha_1-3\alpha_2}
\nonumber\\
&=&1,
\end{eqnarray}
\begin{eqnarray}\label{Ng22}
N_{\alpha_1+\alpha_2,-\alpha_2}&=&-N_{\alpha_1+\alpha_2,-\alpha_1}=-N_{\alpha_1+3\alpha_2,-\alpha_2}
=N_{\alpha_1+3\alpha_2,-\alpha_1-2\alpha_2}=N_{\alpha_1+3\alpha_2,-2\alpha_1-3\alpha_2}
\nonumber\\
&=&-N_{2\alpha_1+3\alpha_2,-\alpha_1}=-N_{2\alpha_1+3\alpha_2,-\alpha_1-\alpha_2}
=N_{2\alpha_1+3\alpha_2,-\alpha_1-2\alpha_2}=N_{2\alpha_1+3\alpha_2,-\alpha_1-3\alpha_2}
\nonumber\\
&=&-N_{-\alpha_1,-\alpha_2}=-N_{-\alpha_1,-\alpha_1-3\alpha_2}=-N_{-\alpha_2,-\alpha_1-2\alpha_2}\nonumber\\
&=&-N_{-\alpha_1-\alpha_2,-\alpha_1-2\alpha_2}=\frac{1}{2},
\end{eqnarray}
\begin{eqnarray}\label{Ng23}
-N_{\alpha_1+2\alpha_2,-\alpha_2}=N_{\alpha_1+2\alpha_2,-\alpha_1-\alpha_2}=-N_{-\alpha_2,-\alpha_1-\alpha_2}=\frac{2}{3}.
\end{eqnarray}
The first relation of (\ref{EEH}) for the root $\alpha_1+\alpha_2$,
as an example, follows from  the last results as
\begin{eqnarray}\label{H21}
H'_{\alpha_1+\alpha_2}&\equiv
&[X'_{\alpha_1+\alpha_2},Y'_{\alpha_1+\alpha_2}]=
[[X'_{\alpha_1},X'_{\alpha_2 }],Y'_{\alpha_1+\alpha_2}]
\nonumber\\
&=&-[[Y'_{\alpha_1+\alpha_2},X'_{\alpha_1}],X'_{\alpha_2}]
-[[X'_{\alpha_2},Y'_{\alpha_1+\alpha_2}],X'_{\alpha_1}]
\nonumber\\
&=&-[Y'_{\alpha_2},X'_{\alpha_2 }]-[Y'_{\alpha_1},X'_{\alpha_1}]
=H'_{\alpha_2}+H'_{\alpha_1}.
\end{eqnarray}
Similarly, we get
\begin{eqnarray}\label{HH}
\begin{array}{ll}
&H^{\prime}_{\alpha_1+2\alpha_2}=H'_{\alpha_1}+2H'_{\alpha_2},\\
&H'_{\alpha_1+3\alpha_2}=H'_{\alpha_1}+3H'_{\alpha_2},\\
&H'_{2\alpha_1+3\alpha_2}=2H'_{\alpha_1}+3H'_{\alpha_2}.
\end{array}
\end{eqnarray}
Again, as an example, for the roots $\alpha=\alpha_1$ and
$\beta=\alpha_2$, the second relation of (\ref{HHEH}) is considered
in the following method
\begin{eqnarray}\label{HXX}
[H'_{\alpha_1},X'_{\alpha_2}]=[[X'_{\alpha_1},Y'_{\alpha_1}],X'_{\alpha_2}]
&=&-[[X'_{\alpha_2},X'_{\alpha_1}],Y'_{\alpha_1}]-[[Y'_{\alpha_1},X'_{\alpha_2}],X'_{\alpha_1}]\nonumber\\
&=&[X'_{\alpha_1+\alpha_2},Y'_{\alpha_1}]=\frac{-1}{2}X'_{\alpha_2}.
\end{eqnarray}
This implies
$\langle\alpha_2,H_{\alpha_1}^{\prime}\rangle=\frac{-1}{2}$ which is
incompatible with (\ref{Cartannumbers}). Therefore, realization of
the constraints (\ref{NN}) and (\ref{NNqp}) violates the negative
integer values for the Cartan numbers. All other structure constants
in the form of $\langle\beta,H_{\alpha}^{\prime}\rangle$ can be
calculated in the same way as in (\ref{HXX}). We recall in Table 2
the unfamiliar commutation relations of the exceptional Lie algebra
${\mathfrak g}_2$ based on the realization of the constraints
(\ref{NN}) and (\ref{NNqp}).

\section{Conclusions}
Finally, we end the paper by comparing Tables 1 and 2. The
commutation relations in Table 1 obey (\ref{N=-N}) and
(\ref{Cartannumbers}) and violate (\ref{NN}) and (\ref{NNqp}). For
Table 2 it is vice versa. In both Tables, the Cartan subalgebra
generators are defined as the first relation of (\ref{EEH}).
Furthermore, linear relations as
\begin{eqnarray}
\begin{array}{llll}
&X'_{\alpha_1}=\sqrt{2}X_{\alpha_1},&
X'_{\alpha_2}=-\frac{1}{\sqrt{3}}X_{\alpha_2},&
X'_{\alpha_1+\alpha_2}=-\sqrt{\frac{2}{3}} X_{\alpha_1+\alpha_2},
\\
&X'_{\alpha_1+2\alpha_2}=2 \frac{\sqrt{2}}{3}X_{\alpha_1+2\alpha_2},
&
X'_{\alpha_1+3\alpha_2}=-2\sqrt{\frac{2}{3}}X_{\alpha_1+3\alpha_2},
&
X'_{2\alpha_1+3\alpha_2}=-\frac{4}{\sqrt{3}}X_{2\alpha_1+3\alpha_2},
\\
&Y'_{\alpha_1}=\frac{1}{2\sqrt{2}}Y_{\alpha_1},&
Y'_{\alpha_2}=\frac{-1}{2\sqrt{3}}Y_{\alpha_2},&
Y'_{\alpha_1+\alpha_2}=-\frac{1}{2\sqrt{6}}Y_{\alpha_1+\alpha_2},
\\
&Y'_{\alpha_1+2\alpha_2}=\frac{1}{4\sqrt{2}}Y_{\alpha_1+2\alpha_2},
&Y'_{\alpha_1+3\alpha_2}=\frac{-\sqrt{6}}{8}Y_{\alpha_1+3\alpha_2},
&Y'_{2\alpha_1+3\alpha_2}=\frac{-\sqrt{3}}{8}Y_{2\alpha_1+3\alpha_2},
\\
&H'_{\alpha_1}=\frac{1}{2}H_{\alpha_1}, &
H'_{\alpha_2}=\frac{1}{6}H_{\alpha_2},&
H'_{\alpha_1+\alpha_2}=\frac{1}{6} H_{\alpha_1+\alpha_2}
\\
&H'_{\alpha_1+2\alpha_2}=\frac{1}{6} H_{\alpha_1+2\alpha_2}, &
H'_{\alpha_1+3\alpha_2}=\frac{1}{2} H_{\alpha_1+3\alpha_2},&
H'_{2\alpha_1+3\alpha_2}=\frac{1}{2} H_{2\alpha_1+3\alpha_2},
\end{array}
\end{eqnarray}
hold an isomorphism between Tabels 1 and 2.

\begin{sidewaystable}
\caption{The commutation relations of the exceptional Lie algebra
$\mathfrak{g}_2$ based on (1.8).} \tiny{
\begin{tabular}{lcccccccccccc}
\toprule
& & & & & & & & & & & &\\
 & $X_{\alpha_1}$ & $X_{\alpha_2}$ & $X_{\alpha_1+\alpha_2}$ & $X_{\alpha_1+2\alpha_2}$
 & $X_{\alpha_1+3\alpha_2}$ & $X_{2\alpha_1+3\alpha_2}$ & $Y_{\alpha_1}$ & $Y_{\alpha_2}$
  & $Y_{\alpha_1+\alpha_2}$ & $Y_{\alpha_1+2\alpha_2}$ & $Y_{\alpha_1+3\alpha_2}$ & $Y_{2\alpha_1+3\alpha_2}$\\
  & & & & & & & & & & & & \\
\toprule
& & & & & & & & & & & &\\
 $H_{\alpha_1}$ & $2X_{\alpha_1}$ & $-X_{\alpha_2}$ &
$X_{\alpha_1+\alpha_2}$ & $0$
 & $-X_{\alpha_1+3\alpha_2}$ & $X_{2\alpha_1+3\alpha_2}$ & $-2Y_{\alpha_1}$ & $Y_{\alpha_2}$
  & $-Y_{\alpha_1+\alpha_2}$ & $0$ & $Y_{\alpha_1+3\alpha_2}$ & $-Y_{2\alpha_1+3\alpha_2}$\\
& & & & & & & & & & & &\\
$H_{\alpha_2}$ & $-3X_{\alpha_1}$ & $2X_{\alpha_2}$ &
$-X_{\alpha_1+\alpha_2}$ & $X_{\alpha_1+2\alpha_2}$
 & $3X_{\alpha_1+3\alpha_2}$ & $0$ & $3Y_{\alpha_1}$ & $-2Y_{\alpha_2}$
  & $Y_{\alpha_1+\alpha_2}$ & $-Y_{\alpha_1+2\alpha_2}$ & $-3Y_{\alpha_1+3\alpha_2}$ & $0$\\
& & & & & & & & & & & &\\
$H_{\alpha_1+\alpha_2}$ & $3X_{\alpha_1}$ & $-X_{\alpha_2}$ &
$2X_{\alpha_1+\alpha_2}$ & $X_{\alpha_1+2\alpha_2}$
 & $0$ & $3X_{2\alpha_1+3\alpha_2}$ & $-3Y_{\alpha_1}$ & $Y_{\alpha_2}$
  & $-2Y_{\alpha_1+\alpha_2}$ & $-Y_{\alpha_1+2\alpha_2}$ & $0$ & $-3Y_{2\alpha_1+3\alpha_2}$\\
& & & & & & & & & & & &\\
$H_{\alpha_1+2\alpha_2}$ & $0$ & $X_{\alpha_2}$ &
$X_{\alpha_1+\alpha_2}$ & $2X_{\alpha_1+2\alpha_2}$
 & $3X_{\alpha_1+3\alpha_2}$ & $-3X_{2\alpha_1+3\alpha_2}$ & $0$ & $-Y_{\alpha_2}$
  & $-Y_{\alpha_1+\alpha_2}$ & $-2Y_{\alpha_1+2\alpha_2}$ & $-3Y_{\alpha_1+3\alpha_2}$ & $-3Y_{2\alpha_1+3\alpha_2}$\\
& & & & & & & & & & & &\\
$H_{\alpha_1+3\alpha_2}$ & $-X_{\alpha_1}$ & $X_{\alpha_2}$ & $0$ &
$X_{\alpha_1+2\alpha_2}$
 & $2X_{\alpha_1+3\alpha_2}$ & $X_{2\alpha_1+3\alpha_2}$ & $Y_{\alpha_1}$ & $-Y_{\alpha_2}$
 & $0$ & $-Y_{\alpha_1+2\alpha_2}$ & $-2Y_{\alpha_1+3\alpha_2}$ & $-Y_{2\alpha_1+3\alpha_2}$\\
& & & & & & & & & & & &\\
$H_{2\alpha_1+3\alpha_2}$ & $X_{\alpha_1}$ & $0$ &
$X_{\alpha_1+\alpha_2}$ & $X_{\alpha_1+2\alpha_2}$
 & $X_{\alpha_1+3\alpha_2}$ & $2X_{2\alpha_1+3\alpha_2}$ & $-Y_{\alpha_1}$ & $0$
 & $-Y_{\alpha_1+\alpha_2}$ & $-Y_{\alpha_1+2\alpha_2}$ & $-Y_{\alpha_1+3\alpha_2}$ & $-2Y_{2\alpha_1+3\alpha_2}$\\
& & & & & & & & & & & &\\
\\ \hline\\
& & & & & & & & & & & &\\
$X_{\alpha_1}$ & $0$ & $X_{\alpha_1+\alpha_2}$ & $0$ & $0$
 & $X_{2\alpha_1+3\alpha_2}$ & $0$ & $H_{\alpha_1}$ & $0$
 & $-Y_{\alpha_2}$ & $0$ & $0$ & $-Y_{\alpha_1+3\alpha_2}$\\
& & & & & & & & & & & &\\
$X_{\alpha_2}$ & $ $ & $0$ & $2X_{\alpha_1+2\alpha_2}$ &
$3X_{\alpha_1+3\alpha_2}$
 & $0$ & $0$ & $0$ & $H_{\alpha_2}$
 & $3Y_{\alpha_1}$ & $-2Y_{\alpha_1+\alpha_2}$ & $-Y_{\alpha_1+2\alpha_2}$ & $0$\\
& & & & & & & & & & & &\\
$X_{\alpha_1+\alpha_2}$ & $ $ & $ $ & $0$ &
$3X_{2\alpha_1+3\alpha_2}$
 & $0$ & $0$ & $-X_{\alpha_2}$ & $3X_{\alpha_1}$
  & $H_{\alpha_1+\alpha_2}$ & $2Y_{2\alpha_2}$ & $0$ & $-Y_{\alpha_1+2\alpha_2}$\\
& & & & & & & & & & & &\\
$X_{\alpha_1+2\alpha_2}$ & $ $ & $ $ & $ $
 & $0$ & $0$ & $0$ & $0$
 & $-2X_{\alpha_1+\alpha_2}$ & $2X_{\alpha_2}$ & $H_{\alpha_1+2\alpha_2}$ & $Y_{\alpha_2}$
& $Y_{\alpha_1+\alpha_2}$\\
& & & & & & & & & & & &\\
$X_{\alpha_1+3\alpha_2}$ & $ $ & $ $ & $ $
 & $ $ & $0$ & $0$ & $0$
 & $-X_{\alpha_1+2\alpha_2}$ & $0$ & $X_{\alpha_2}$ & $H_{\alpha_1+3\alpha_2}$
& $Y_{\alpha_1}$\\
& & & & & & & & & & & &\\
$X_{2\alpha_1+3\alpha_2}$ & $ $ & $ $ & $ $
 & $ $ & $ $ & $0$ & $-X_{\alpha_1+3\alpha_2}$
 & $0$ & $-X_{\alpha_1+2\alpha_2}$ & $X_{\alpha_1+\alpha_2}$ & $X_{\alpha_1}$
& $H_{2\alpha_1+3\alpha_2}$\\
& & & & & & & & & & & &\\
\\ \hline\\
& & & & & & & & & & & &\\
$Y_{\alpha_1}$ & $ $ & $ $ & $ $
 & $ $ & $ $ & $ $ & $0$
 & $-Y_{\alpha_1+\alpha_2}$ & $0$ & $0$ & $-Y_{2\alpha_1+3\alpha_2}$
& $0$\\
& & & & & & & & & & & &\\
$Y_{\alpha_2}$ & $ $ & $ $ & $ $
 & $ $ & $ $ & $ $ & $ $
 & $0$ & $-2Y_{\alpha_1+\alpha_2}$ & $-3Y_{\alpha_1+3\alpha_2}$ & $0$
& $0$\\
& & & & & & & & & & & &\\
$Y_{\alpha_1+\alpha_2}$ & $ $ & $ $ & $ $
 & $ $ & $ $ & $ $ & $ $
 & $ $ & $0$ & $-3Y_{2\alpha_1+3\alpha_2}$ & $0$
& $0$\\
& & & & & & & & & & & &\\
$Y_{\alpha_1+2\alpha_2}$ & $ $ & $ $ & $ $
 & $ $ & $ $ & $ $ & $ $
 & $ $ & $ $ & $0$ & $0$ & $0$\\
& & & & & & & & & & & &\\
$Y_{\alpha_1+3\alpha_2}$ & $ $ & $ $ & $ $
 & $ $ & $ $ & $ $ & $ $
 & $ $ & $ $ & $ $ & $0$ & $0$\\
& & & & & & & & & & & &\\
$Y_{2\alpha_1+3\alpha_2}$ &  &  &
 &  &  &  &
 &  &  &  &  & $0$\\
 & & & & & & & & & & & &\\
\hline\\
\end{tabular}
}
\end{sidewaystable}

\begin{sidewaystable}
\caption{The commutation relations of the exceptional Lie algebra
${\mathfrak g}_2$ based on (\ref{NN}) and (\ref{NNqp}). }\tiny{
\begin{tabular}{lcccccccccccc}
\toprule
& & & & & & & & & & & &\\
$ $ & $X'_{\alpha_1}$ & $X'_{\alpha_2}$ & $X'_{\alpha_1+\alpha_2}$ &
$X'_{\alpha_1+2\alpha_2}$ & $X'_{\alpha_1+3\alpha_2}$ &
$X'_{2\alpha_1+3\alpha_2}$ & $Y'_{\alpha_1}$ & $Y'_{\alpha_2}$
& $Y'_{\alpha_1+\alpha_2}$ & $Y'_{\alpha_1+2\alpha_2}$ & $Y'_{\alpha_1+3\alpha_2}$ & $Y'_{2\alpha_1+3\alpha_2}$\\
& & & & & & & & & & & &\\
\toprule
& & & & & & & & & & & &\\
 $H'_{\alpha_1}$ & $X'_{\alpha_1}$ &
$\frac{-1}{2}X'_{\alpha_2}$ & $\frac{1}{2}X'_{\alpha_1+\alpha_2}$ &
$0$ & $\frac{-1}{2}X'_{\alpha_1+3\alpha_2}$ &
$\frac{1}{2}X'_{2\alpha_1+3\alpha_2}$ & $-Y'_{\alpha_1}$ &
$\frac{-1}{2}Y'_{\alpha_2}$
& $\frac{-1}{2}Y'_{\alpha_1+\alpha_2}$ & $0$ & $\frac{1}{2}Y'_{\alpha_1+3\alpha_2}$ & $\frac{-1}{2}Y'_{2\alpha_1+3\alpha_2}$\\
& & & & & & & & & & & &\\
$H'_{\alpha_2}$ & $\frac{-1}{2}X'_{\alpha_1}$ &
$\frac{1}{3}X'_{\alpha_2}$ & $\frac{-1}{6}X'_{\alpha_1+\alpha_2}$ &
$\frac{1}{6}X'_{\alpha_1+2\alpha_2}$ &
$\frac{1}{2}X'_{\alpha_1+3\alpha_2}$ & $0$ &
$\frac{1}{2}Y'_{\alpha_1}$ & $\frac{-1}{3}Y'_{\alpha_2}$
& $\frac{1}{6}Y'_{\alpha_1+\alpha_2}$ & $\frac{-1}{6}Y'_{\alpha_1+2\alpha_2}$ & $\frac{-1}{2}Y'_{\alpha_1+3\alpha_2}$ & $0$\\
& & & & & & & & & & & &\\
$H'_{\alpha_1+\alpha_2}$ & $\frac{1}{2}X'_{\alpha_1}$ &
$\frac{-1}{6}X'_{\alpha_2}$ & $\frac{1}{3}X'_{\alpha_1+\alpha_2}$ &
$\frac{1}{6}X'_{\alpha_1+2\alpha_2}$ & $0$ &
$\frac{1}{2}X'_{2\alpha_1+3\alpha_2}$ & $\frac{-1}{2}Y'_{\alpha_1}$
& $\frac{1}{6}Y'_{\alpha_2}$ & $\frac{-1}{3}Y'_{\alpha_1+\alpha_2}$
& $\frac{-1}{6}Y'_{\alpha_1+2\alpha_2}$
& $0$ & $\frac{-1}{2}Y'_{2\alpha_1+3\alpha_2}$\\
& & & & & & & & & & & &\\
$H'_{\alpha_1+2\alpha_2}$ & $0$ & $\frac{1}{6}X'_{\alpha_2}$ &
$\frac{1}{6}X'_{\alpha_1+\alpha_2}$ &
$\frac{1}{3}X'_{\alpha_1+2\alpha_2}$ &
$\frac{1}{2}X'_{\alpha_1+3\alpha_2}$ &
$\frac{1}{2}X'_{2\alpha_1+3\alpha_2}$ & $0$ &
$\frac{-1}{6}Y'_{\alpha_2}$ & $\frac{-1}{6}Y'_{\alpha_1+\alpha_2}$ &
$\frac{-1}{3}Y'_{\alpha_1+2\alpha_2}$
& $\frac{-1}{2}Y'_{\alpha_1+3\alpha_2}$ & $\frac{-1}{2}Y'_{2\alpha_1+3\alpha_2}$\\
& & & & & & & & & & & &\\
$H'_{\alpha_1+3\alpha_2}$ & $\frac{-1}{2}X'_{\alpha_1}$ &
$\frac{1}{2}X'_{\alpha_2}$ & $0$ &
$\frac{1}{2}X'_{\alpha_1+2\alpha_2}$ & $X'_{\alpha_1+3\alpha_2}$ &
$\frac{1}{2}X'_{2\alpha_1+3\alpha_2}$ & $\frac{1}{2}Y'_{\alpha_1}$ &
$\frac{-1}{2}Y'_{\alpha_2}$ & $0$ &
$\frac{-1}{2}Y'_{\alpha_1+2\alpha_2}$ & $-Y'_{\alpha_1+3\alpha_2}$
& $\frac{-1}{2}Y'_{2\alpha_1+3\alpha_2}$\\
& & & & & & & & & & & &\\
$H'_{2\alpha_1+3\alpha_2}$ & $\frac{1}{2}X'_{\alpha_1}$ & $0$ &
$\frac{1}{2}X'_{\alpha_1+\alpha_2}$ &
$\frac{1}{2}X'_{\alpha_1+2\alpha_2}$ &
$\frac{1}{2}X'_{\alpha_1+3\alpha_2}$ & $X'_{2\alpha_1+3\alpha_2}$ &
$\frac{-1}{2}Y'_{\alpha_1}$ & $0$ &
$\frac{-1}{2}Y'_{\alpha_1+\alpha_2}$ &
$\frac{-1}{2}Y'_{\alpha_1+2\alpha_2}$
& $\frac{-1}{2}Y'_{\alpha_1+3\alpha_2}$ & $-Y'_{2\alpha_1+3\alpha_2}$\\
& & & & & & & & & & & &\\
\\ \hline\\
& & & & & & & & & & & &\\
$X'_{\alpha_1}$ & $0$ & $X'_{\alpha_1+\alpha_2}$ & $0$ & $0$ &
$X'_{2\alpha_1+3\alpha_2}$ & $0$ & $H'_{\alpha_1}$ & $0$
& $-Y'_{\alpha_2}$ & $0$ & $0$ & $-Y'_{\alpha_1+3\alpha_2}$\\
& & & & & & & & & & & &\\
$X'_{\alpha_2}$ & & $0$ & $X'_{\alpha_1+2\alpha_2}$ &
$X'_{\alpha_1+3\alpha_2}$ & $0$ & $0$ & $0$ & $H'_{\alpha_2}$
& $Y'_{\alpha_1}$ & $-Y'_{\alpha_1+\alpha_2}$ & $-Y'_{\alpha_1+2\alpha_2}$ & $0$\\
& & & & & & & & & & & &\\
$X'_{\alpha_1+\alpha_2}$ & & & $0$ & $X'_{2\alpha_1+3\alpha_2}$ &
$0$ & $0$ & $\frac{-1}{2}X'_{\alpha_2}$ & $\frac{1}{2}X'_{\alpha_1}$
& $H'_{\alpha_1+\alpha_2}$ & $Y'_{\alpha_2}$ & $0$ & $-Y'_{\alpha_1+2\alpha_2}$\\
& & & & & & & & & & & &\\
$X'_{\alpha_1+2\alpha_2}$ & & & & $0$ & $0$ & $0$ & $0$ &
$\frac{-2}{3}X'_{\alpha_1+\alpha_2}$ & $\frac{2}{3}X'_{\alpha_2}$ &
$H'_{\alpha_1+2\alpha_2}$ & $Y'_{\alpha_2}$
& $Y'_{\alpha_1+\alpha_2}$\\
& & & & & & & & & & & &\\
$X'_{\alpha_1+3\alpha_2}$ & & & & & $0$ & $0$ & $0$ &
$\frac{-1}{2}X'_{\alpha_1+2\alpha_2}$ & $0$ &
$\frac{1}{2}X'_{\alpha_2}$ & $H'_{\alpha_1+3\alpha_2}$
& $\frac{1}{2}Y'_{\alpha_1}$\\
& & & & & & & & & & & &\\
$X'_{2\alpha_1+3\alpha_2}$ & & & & & & $0$ &
$\frac{-1}{2}X'_{\alpha_1+3\alpha_2}$ & $0$ &
$\frac{-1}{2}X'_{\alpha_1+2\alpha_2}$ &
$\frac{1}{2}X'_{\alpha_1+\alpha_2}$ & $\frac{1}{2}X'_{\alpha_1}$
& $H'_{2\alpha_1+3\alpha_2}$\\
& & & & & & & & & & & &\\
\\ \hline\\
& & & & & & & & & & & &\\
$Y'_{\alpha_1}$ & & & & & & & $0$
& $\frac{-1}{2}Y'_{\alpha_1+\alpha_2}$ & $0$ & $0$ & $\frac{-1}{2}Y'_{2\alpha_1+3\alpha_2}$ & $0$\\
& & & & & & & & & & & &\\
$Y'_{\alpha_2}$ & & & & & & &
& $0$ & $\frac{-2}{3}Y'_{\alpha_1+\alpha_2}$ & $\frac{-1}{2}Y'_{\alpha_1+3\alpha_2}$ & $0$ & $0$\\
& & & & & & & & & & & &\\
$Y'_{\alpha_1+\alpha_2}$ & & & & & & & & & $0$ & $\frac{-1}{2}Y'_{2\alpha_1+3\alpha_2}$ & $0$ & $0$\\
& & & & & & & & & & & &\\
$Y'_{\alpha_1+2\alpha_2}$ & & & & & & & & & & $0$ & $0$ & $0$\\
& & & & & & & & & & & &\\
$Y'_{\alpha_1+3\alpha_2}$ & & & & & & & & & & & $0$ & $0$\\
& & & & & & & & & & & &\\
$Y'_{2\alpha_1+3\alpha_2}$ & & & & & & & & & & & & $0$\\
& & & & & & & & & & & &\\
\hline\\
\end{tabular}
}
\end{sidewaystable}

\end{document}